\documentclass[a4paper]{amsart}

\usepackage{amssymb,amsmath}
\usepackage{latexsym}
\usepackage{eucal}
\makeatletter
\@addtoreset{equation}{section}\makeatother

\newtheorem{theo}{Theorem}[section]
\newtheorem{lem}[theo]{Lemma}
\newtheorem{prop}[theo]{Proposition}

\newcommand{\mc}{\mathcal}
\newcommand{\rr}{\mathbb{R}}
\newcommand{\nn}{\mathbb{N}}
\newcommand{\cc}{\mathbb{C}}
\newcommand{\hh}{\mathbb{H}}
\newcommand{\zz}{\mathbb{Z}}

\newcommand{\la}{\lambda}
\newcommand{\eps}{\epsilon}

\newcommand{\pl}{\partial}
\newcommand{\x}{\times}

\newcommand{\til}{\widetilde}

\newcommand{\cjd}{\rangle}
\newcommand{\cjg}{\langle}

\newcommand{\demi}{\frac{1}{2}}
\newcommand{\ndemi}{\frac{n}{2}}
\newcommand{\tra}{\textrm{Tr}}

\newcommand{\rang}{\textrm{rank}}

\newcommand{\indic}{\operatorname{1\negthinspace l}}
\newcommand{\zerotr}{\mathop{\hbox{\rm 0-tr}}\nolimits}
\newcommand{\zerov}{\mathop{\hbox{\rm 0-vol}}\nolimits}
\def\qed{\hfill$\square$}
\begin{document}
\title[Wave $0$-trace and length spectrum]
{Wave $0$-trace and length spectrum on convex co-compact hyperbolic manifolds}

\author[Colin Guillarmou]{Colin Guillarmou}
\address{Department of mathematics\\
         Purdue University\\
         150 N. University Street, West-Lafayette\\
         IN 47907, USA}
     \email{cguillar@math.purdue.edu}
\author[Fr\'ed\'eric Naud]{Fr\'ed\'eric Naud}
\address{Department of Mathematics\\
University of California Berkeley\\
Berkeley CA 94720-3840 USA}
\email{fnaud@math.berkeley.edu}

%\thanks{}
%\keywords{}
%
\begin{abstract}
\noindent For convex co-compact hyperbolic quotients $\Gamma\backslash\hh^{n+1}$
we obtain a formula relating the $0$-trace of
the wave operator with the resonances and some conformal invariants of the boundary,
generalizing a formula of Guillop\'e and Zworski in dimension $2$.
Then, by writing this $0$-trace with the length spectrum, we prove precise asymptotics of the
number of closed geodesics with an effective, exponentially small error term
when the dimension of the limit set of $\Gamma$ is greater than $\ndemi$.
\end{abstract}
\maketitle

%\newpage
%\tableofcontents
%\newpage

\section{Introduction}

The purpose of this note is to establish the relation between the
renormalized trace (called $0$-trace) of the wave operator,
the resonances, some conformal invariants of the boundary and
the length spectrum on a
convex co-compact quotient of the hyperbolic space $\hh^{n+1}$.
From the work of Patterson-Perry \cite{PP} and Bunke-Olbrich \cite{BO} on Selberg zeta function, we derive
a trace formula similar to the one given by Guillop\'e and Zworski \cite{GZ,GZ1} on surfaces.
Then we compute the $0$-trace of the wave operator in terms of the primitive geodesics
to obtain asymptotic expansions for the counting function of the set of closed geodesics with precise
error terms related to the spectrum of the Laplacian, improving
asymptotics previously obtained by Perry \cite{P4}.\\

Actually, a Selberg formula relating resonances and length spectrum
has already been obtained by Perry \cite{P} in this setting but Perry did not
use explicitly the $0$-trace of the wave operator. The first motivation for studying
this $0$-trace is that our formula could be extended to more general
settings like asymptotically hyperbolic manifolds.
Another interest of this formula comes from the conformal invariants of the boundary
which appear in the expansion of the $0$-trace of the wave operator
and are related to the divisors of Selberg's zeta function at certain values.\\

Let us first recall standard notations and definitions. A discrete
group $\Gamma$ of (orientation preserving) isometries of the $n+1$-dimensional real hyperbolic
space $\hh^{n+1}$ is called convex co-compact if it admits a finite
sided fundamental domain whose intersection with $\partial \hh^{n+1}$
does not touch the limit set of $\Gamma$. If we require that $\Gamma$
has no elliptic elements, then $X=\Gamma \backslash \hh^{n+1}$ is a
hyperbolic manifold of infinite volume. The non-wandering set of the
geodesic flow on the unit tangent bundle is then a compact set whose
Hausdorff dimension is $2\delta+1$, where $\delta$ is
the dimension of the limit set, see \cite{ZW}. The dimension $\delta$ is also the
topological entropy of the geodesic flow on its non-wandering set, see
the work of Sullivan \cite{S1,S2}. Let $\mc{P}$ be the set of
primitive closed geodesics $\gamma$ on $X$, and if $\gamma \in
\mc{P}$, let $l(\gamma)$ denote its length.\\

A wide class of convex co-compact groups are given by Schottky
groups of isometries, and every convex co-compact surface is in fact
obtained by a quotient of $\hh^2$ by a fuchsian Schottky group.
In higher dimensions, as pointed out by
Maskit \cite{Mask}, the set of Schottky groups does not exhaust the set of
convex co-compact groups. Indeed, a 3-manifold is Schottky if and only if its fundamental group
is a free group of finite type, thus all convex co-compact manifolds obtained by a quotient of $\hh^3$ by
quasifuchsian groups cannot be Schottky.

We also recall that on a convex co-compact hyperbolic quotient $X=\Gamma\backslash\hh^{n+1}$
equipped with its hyperbolic metric $g$,
the Laplacian $\Delta_X$ has for continuous spectrum the half line $[\frac{n^2}{4},+\infty)$
and a finite set $\sigma_{pp}(\Delta_X)$ of eigenvalues included in
$(0,\frac{n^2}{4})$.
The modified resolvent of the Laplacian
\[R_X(s):=(\Delta_X-s(n-s))^{-1}\]
is  meromorphic on $\{\Re(s)>\frac{n}{2}\}$ with a finite number of
poles related to $\sigma_{pp}(\Delta_X)$,
and extends to $s\in\cc$ (as a map from $L^2_{comp}(X)$ to $L^2_{loc}(X)$) with poles
of finite multiplicity called resonances (see \cite{MM,GZ2}), the multiplicity of a resonance $s_0$ being defined by
\begin{equation}\label{multip1}
m_{s_0}:=\rang (\textrm{Res}_{s_0}R_X(s)).
\end{equation}
We will denote by 
\[\mc{R}:=\{s_0\in\cc; R(s) \textrm{ has a pole at }s=s_0 \} \]
the set of resonances of $\Delta_X$.
We can also add that $X$ can be compactified in a compact manifold
with boundary $\bar{X}$ such that, for $x$ a boundary defining
function of $\pl\bar{X}$ in $\bar{X}$ (i.e. $\pl\bar{X}=\{x=0\}$
and $dx|_{\pl\bar{X}}\not=0$), the metric $x^2g$ on $X$ extends
smoothly to $\bar{X}$. In view of the non-uniqueness of the
boundary defining function, the boundary carries a natural
conformal class of metric $[h_0]$ associated to $g$ by taking the
conformal class of $h_0:=x^2g|_{\pl\bar{X}}$. Graham, Jenne,
Manson and Sparling \cite{GJMS} have defined some natural
conformally covariant powers of the Laplacian on conformal
manifolds $(M,[h_0])$ and those were identified by Graham and
Zworski \cite{GRZ} with poles of the scattering operator on
asymptotically hyperbolic Einstein manifolds with conformal
infinity $(M,[h_0])$ by using scattering theory. We will denote by
$P_k$ the k-th conformal power of the Laplacian on the conformal
infinity $(\pl\bar{X},[h_0])$. Note that in even dimension of
$\pl\bar{X}$, the existence of natural $P_k$ for $k>\ndemi$ does
not hold (see Gover-Hirachi \cite{GOH}) and the definition of what
we call $P_k$ will be detailed later in that case.

On a convex co-compact hyperbolic manifold,
although the trace of the wave operator is not well defined,
there is a natural renormalization called $0$-trace. Thus,
following Joshi-Sa Barreto \cite{JSB1} or  Guillop\'e-Zworski \cite{GZ,GZ1},
we can define the $0$-trace of the wave operator as a distribution on $\rr^*$ noted
\[\zerotr\left(\cos\left(t\sqrt{\Delta_{X}-\frac{n^2}{4}}\right)\right).\]
The definition of this renormalized trace will be detailed in the next section.

We recall that the set of closed geodesics $\mc{P}$ is in one-to-one
correspondence with primitive conjugacy classes of hyperbolic
isometries in $\Gamma$. Given an hyperbolic element $h\in \Gamma$
associated to a closed geodesic $\gamma$,
then there exists $\alpha \in {\rm Isom}(\hh^{n+1})$ such that for all
$(x,y)\in \hh^{n+1}=\rr^n\times \rr_+$,
$$\alpha^{-1}\circ h \circ \alpha(x,y)=e^{l(\gamma)}(O_\gamma(x),y),$$
where $O_\gamma \in SO_n(\rr)$. We will denote by
$\alpha_1(\gamma),\ldots,\alpha_n(\gamma)$ the eigenvalues of
$O_\gamma$, and we set
\footnote{If $P_\gamma$ is the Poincar\'e linear map asociated to the
primitive periodic orbit $\gamma$ of the geodesic flow on the unit
tangent bundle, then $$|\det\left( I-P_\gamma^m
\right)|^{1/2}=e^{\ndemi m l(\gamma)}G_\gamma(m).$$}
$$G_\gamma(k)=\det \left(I-e^{-kl(\gamma)}O_\gamma^k \right)=
\prod_{i=1}^n \left(1-e^{-kl(\gamma)}\alpha_i(\gamma)^k \right).$$

We can now state the wave $0$-trace formula.

\begin{theo}\label{poisson}
Let $X=\Gamma\backslash\hh^{n+1}$ be a convex co-compact hyperbolic manifold, $(P_k)_{k\in\nn}$
the k-th conformal Laplacian on the boundary $\pl\bar{X}$, then
we have as distributions on $\rr^*$
\[\zerotr\left(\cos\left(t\sqrt{\Delta_{X}-\frac{n^2}{4}}\right)\right)=
\demi\sum_{s\in\mc{R}}m_s e^{-(\ndemi-s)|t|}+\demi\sum_{k\in\nn}d_ke^{-k|t|}
-2^{-n-1}A(X)\frac{\cosh\frac{t}{2}}{(\sinh\frac{|t|}{2})^{n+1}},\]
where
\[d_k:=\dim\ker P_k,\quad
A(X)=\left\{\begin{array}{ll}
0 & \textrm{ if }n+1\textrm{ is even}\\
\chi(X)\textrm{ the Euler characteristic of }X & \textrm{ if } n+1\textrm{ is odd}\end{array}\right..\]
We also have
\[\zerotr\left(\cos\left(t\sqrt{\Delta_{X}-\frac{n^2}{4}}\right)\right)=
\sum_{\gamma}\sum_{m=1}^{\infty}\frac{l(\gamma)e^{-\ndemi ml(\gamma)}
}{2G_\gamma(m)}\delta(|t|-ml(\gamma))+B(X)\frac{\cosh\frac{t}{2}}
{(\sinh\frac{|t|}{2})^{n+1}}\]
where the sum runs over the primitive closed geodesics of $X$ and
\[B(X)=\left\{\begin{array}{ll}
n!!2^{-\frac{3(n+1)}{2}}(-\pi)^{-\frac{n+1}{2}}\zerov(X) & \textrm{ if }n+1\textrm{ is even}\\
0 & \textrm{ if } n+1\textrm{ is odd.}\end{array}\right.\]
\end{theo}

\textsl{Remark}: the first formula is different from Euclidean cases in view of these additional
terms $d_k$ coming from the boundary.
Of course, we expect that this one could be extended to compact perturbations of hyperbolic convex
co-compact manifolds and possibly to asymptotically Einstein manifolds,
at least when $n+1$ is even.
Observe also that it  corresponds to Guillop\'e-Zworski trace formula in \cite{GZ} in the sense that
in dimension $2$ we have $d_k=0$ for all $k$ (see \cite[Prop. 4.3]{BO1} where in this case $d_k$ must be interpreted as
the dimension of the kernel of the residue of the scattering matrix at $\ndemi+k$).\\

As a corollary of the trace formula, we can study the asymptotic behaviour of the primitive geodesics
counting function. We recall that $\mc{P}$ denotes the set of primitive closed geodesics, the set 
$\{l(\gamma); \gamma\in \mc{P}\}$ is often refered as the length spectrum of the manifold 
$X=\Gamma\backslash\hh^{n+1}$, it is a countable subset of $\rr^+$ which
accumulates only at $+\infty$ (for more details about the length spectrum of hyperbolic manifolds,
we refer the reader to Buser \cite{Bu} or Hejhal \cite{Hej}).
The basic counting function $N(T)$ for the closed geodesics is defined
as usual for $T \geq 0$ by
$$N(T):=\# \{\gamma \in \mc{P}\ ;\ l(\gamma)\leq T \},$$

In this paper, we show the following.
\begin{theo}
\label{lsp}
Let $X=\Gamma\backslash\hh^{n+1}$ be a convex co-compact manifold as
above such that $\delta >\frac{n}{2}$. As $T\rightarrow
+\infty$, we have
$$N(T)=\rm{li}(e^{\delta T})+\sum_{\beta_n(\delta)<\alpha_i<\delta}
\rm{li}(e^{\alpha_i T})+O\left( \frac{e^{\beta_n(\delta)T}}{T}\right),$$
where $\rm{li}(x)=\int_2^x\frac{dt}{\log(t)}$ and
$\beta_n(\delta)=\frac{n}{n+1}(\frac{1}{2}+\delta)$.
The coefficients $\alpha_i$ are in bijection with
$\sigma_{pp}(\Delta_X)$ by $\alpha_i(n-\alpha_i)=\lambda_i\in
\sigma_{pp}(\Delta_X)$.
\end{theo}
\textsl{Remarks}. The leading term $N(T)\sim \frac{e^{\delta T}}{\delta T}$
has been obtained by Perry in \cite{P4} without any assumption on the
dimension $\delta$. If we let $\delta\rightarrow 1$ then we
essentially recover the error term $O(e^{\frac{3}{4}T}/T)$ known for
finite area surfaces, see the work of Huber, Hejhal, Randol and Sarnak \cite{Hub, Hej, Rand, Sar}.
For 3-manifolds, Theorem \ref{lsp} can be applied for example to the case of strictly
quasifuchsian groups $\Gamma$ (i.e. the limit set is a strict quasicircle)
for which we know since Bowen \cite{Bo} that $\delta>1=\ndemi$.
In the particular case of surfaces ($n=1$), the second author \cite{N1} has obtained unconditional
exponentially small error terms using transfer operator techniques and
thermodynamical formalism. We point out that this asymptotic expansion also holds
for geometrically finite surfaces with cusps \cite{N2}.
Very recently, a non-trivial extension of these
techniques has been obtained by L. Stoyanov \cite{Stoy}. His result
implies an exponentially small error term for the counting function $N(T)$
for an arbitrary convex co-compact manifold without any assumptions on
$\delta$. However, because of the techniques employed, his error term is
not explicit and has no spectral interpretation.\\

\textbf{Acknowledgements:} the first author would like to thank Robin Graham for helpful
discussions and Siu-Hung Tang for pointing out the appendix of Patterson-Perry \cite{PP}
to compute the conformal Laplacians on hyperbolic compact manifolds. Also, we thank
Laurent Guillop\'e for comments on this work. First author acknowledges support 
of NSF grant DMS0500788.

\section{0-trace formula}

\subsection{0-Renormalization}
Let $X=\Gamma\backslash \hh^{n+1}$ a convex co-compact hyperbolic
manifold equipped with the hyperbolic metric $g$ and let $\bar{X}$
be its compactification as a smooth manifold with boundary. On
such a manifold and after having chosen a boundary defining
function $x$, we can define the 0-integral of a smooth function
$f$ on $\bar{X}$ by the formula
\[\int^0_X f:=FP_{\eps\to 0}\int_{x(m)>\eps}f(m)dvol_{g}(m),\]
this depends of course on the function $x$ and this definition can be extended
for functions for which the finite part exists.
However, it is shown for exemple by Graham \cite{GR}
that the $0$-volume of $X$, defined as the $0$-integral of the function $1$,
is independent of the choice of $x$ (here $X$ is Einstein) if the dimension $n+1$ of $X$ is
even. For an operator $A$ on $X$ which has a Schwartz kernel $A(m,m')$ which
is smooth when restricted to the diagonal of $X\x X$, one can also define its $0$-trace by
\[\zerotr A := \int^0_X A(w,w)\]
when the $0$-integral exists (see \cite{GZ,JSB1,PP,BJP} for details and examples).
Note that those renormalizations are naturally related to the $0$-calculus and $0$-structure
defined by Mazzeo-Melrose \cite{MM,MA1,MA2} (i.e. related to the
`geometric operators' on conformally compact manifolds).

Following Joshi-Sa Barreto \cite{JSB1}
(see also Guillop\'e-Zworski \cite{GZ,GZ1}),
we can define its $0$-trace as a distribution on $\rr^*$
\begin{equation}\label{otrace}
u(t):=\zerotr\left(\cos\left(t\sqrt{\Delta_{X}-\frac{n^2}{4}}\right)\right).
\end{equation}
which means that for all $\varphi\in C_0^\infty(\rr^*)$
\[\cjg u,\varphi\cjd:= \zerotr \cjg \cos\left(\bullet\sqrt{\Delta_{X}-\frac{n^2}{4}}\right),\varphi\cjd.\]

\subsection{Resonances, scattering poles and conformal operators}
If $s$ is not a resonance and not in $\demi\zz$,
the generalized eigenfunctions of $\Delta_X$ for the `eigenvalue' $s(n-s)$
have the following behaviour on $\bar{X}$
\[E(s)=x^{s}F_1(s)+x^{n-s}F_2(s), \quad F_i(s)\in C^{\infty}(\bar{X})\]
and one can show (see \cite{GRZ} for example) that it is unique if one requires that
$F_2(s)|_{\pl\bar{X}}=f_0$ for a fixed function $f_0\in C^\infty(\pl\bar{X})$.
Thus one can define the scattering operator as the operator on $C^\infty (\pl\bar{X})$
\[S(s): f_0\to F_1(s)|_{\pl\bar{X}}.\]
It turns out that (see \cite{P1,PP}) for $h_0:=(x^2g)|_{\pl\bar{X}}$ the operator
\[\til{S}(s):=2^{2s-n}\frac{\Gamma(s-\ndemi)}{\Gamma(\ndemi-s)}(1+\Delta_{h_0})^{\frac{-s+n}{2}}
S(s)(1+\Delta_{h_0})^{\frac{-s+n}{2}}\]
has also a meromorphic extension to $\cc$ with poles (called scattering poles) of finite multiplicity
in $\cc$, the multiplicity of a pole $s_0$ being defined here by
\begin{equation}\label{multip2}
\nu_{s_0}:=-\tra(\textrm{Res}_{s_0}(\til{S}'(s)\til{S}^{-1}(s))).
\end{equation}
The set of the scattering poles contained in $\{\Re(s)<\ndemi\}$
will be denoted by $\mc{S}$ whereas the set of resonances will be
denoted by $\mc{R}$. In \cite{G}, we have given a formula relating
the multiplicities (\ref{multip1}) and (\ref{multip2}) on an
general asymptotically hyperbolic manifold. We recall this result
in the present case and extend it slightly (see the following
paragraph for the precise definition of what exactly is $P_k$ in
this setting):
\begin{prop}\label{multip}
Let $X=\Gamma\backslash \hh^{n+1}$ be a convex co-compact hyperbolic manifold with conformal
infinity $(\pl\bar{X},[h_0])$ and let $\Re(s_0)<\ndemi$, then
we have the relation
\begin{equation}\label{multi}
\nu_{s_0}=m_{s_0}-m_{n-s_0}+\indic_{\ndemi-\nn}(s_0)\dim\ker P_{\ndemi-s_0}
\end{equation}
with $\indic_{\ndemi-\nn}$ the characteristic function of $\ndemi-\nn$ and $P_k$ for $k\in\nn$
is the k-th conformal Laplacian on the conformal manifold $(\pl\bar{X},[h_0])$ defined in
Graham-Zworski \cite{GRZ}.
\end{prop}
\textsl{Proof}: in \cite{G} we dealt with all cases except when
$s_0\in\{s; s(n-s)\in\sigma_{pp}(\Delta_g)\}\cap (\ndemi-\nn)$. In fact we can deal with these special points
using the perturbation method of Borthwick-Perry \cite{PB}. Indeed, since the resolvent
and scattering operator are meromorphic with poles of finite multiplicity
in $\{\Re(s)<\ndemi\}$, one can remark from \cite{PB} that if $s_0\in\frac{n-k}{2}$, it is possible to add
a sufficiently small non-negative compactly supported potential $V$ on $X$, such that $m_{n-s_0}$, $m_{s_0}$ and $\nu_{s_0}$
remains invariant and the eigenvalue $s_0(n-s_0)$ is pushed a little so that
$s_0(n-s_0)\notin\sigma_{pp}(\Delta_g+V)$.
The formula being now satisfied at $s_0$ we have the result since $P_\frac{k}{2}$ when $k$ is even
depends only on the $k$ first derivatives $(\pl^j_x(x^2g)|_{x=0})_{j=1\dots,k}$ if $g$ is the hyperbolic
metric on $X$ (see again \cite{GRZ}).
\qed\\

Now we say a few words about the conformal powers of the
Laplacian on the compact conformal manifold $(\pl\bar{X},[h_0])$.
By \cite{GJMS}, for each conformal representative $h$ in a
conformal class $[h_0]$ on $\pl\bar{X}$, the $k$-th conformal
Laplacian $P_k$ associated to $h$ is, for all
$k\in\nn$ if $n$ is odd (resp. for $k\leq \ndemi$ if $n$ is even), a well defined 
natural differential operator of order $2k$ with principal symbol
\[\sigma_0(P_k)=\sigma_0(\Delta^k_{h})\]
which satisfies a covariant rule when the conformal representative
is changed, that is
\begin{equation}\label{covariant}
\hat{P}_k=e^{-(\ndemi+k)\omega}P_ke^{(\ndemi-k)\omega}
\end{equation}
if $\hat{h}=e^{2\omega}h$ (for $\omega\in C^\infty(\pl\bar{X})$)
and $\hat{P}_k$ the operator associated to $\hat{h}$. The first
example is the conformal Laplacian
$P_1=\Delta_{h}+\frac{n-2}{4(n-1)}R$ where $R$ is the scalar
curvature of $h$ on $\pl\bar{X}$ and $P_2$ is the so-called
Paneitz operator. If $n=1$ or if $n$ is even and $k>\ndemi$ these
conformal Laplacians $P_k$ are not naturally well-defined in
general but for locally conformally flat manifolds, it is also
proved in \cite{GJMS} that the $k$-th conformal Laplacian is
well-defined canonically for all $k$ and $n>2$, including when $n$
is even. With the restriction $k\leq\ndemi$ if $n$ is even, Graham
and Zworski proved that $c_kP_k$ is equal to $c_kp_k$ with
$c_k=(-1)^{k+1}(2^{2k}k!(k-1)!)^{-1}$ and $p_k$ the residue of the
scattering operator at $\ndemi+k$ (modulo smoothing operators on
$\pl\bar{X}$ involving only the possible $L^2$ eigenvalues of
$\Delta_X$), where we recall that we are in the framework of
asymptotically hyperbolic Einstein manifolds. Actually $X$ is
hyperbolic, thus the conformal infinity $(\pl\bar{X},[h_0])$ is
locally conformally flat and, according to Robin Graham
\cite{GRpr}, the operator $c_k^{-1}p_k$ is the canonical $P_k$
constructed in \cite{GJMS} on locally conformally flat manifolds
for $n>2$ even and all $k\in\nn$. In any case and by convention,
$P_k$ will be called $k$-th conformal power of the Laplacian and
defined as the differential operator $P_k:=c_k^{-1}p_k$ and
$d_k=:\dim\ker P_k$ is a finite number which is conformally
invariant on the conformal infinity since $P_k$ is Fredholm and
satisfies a covariant rule.

There is a special case where the operators $P_k$ are computable,
this is when the conformal infinity has a conformal representative
with constant curvatures $K$. Note that by dilation of the metric
and the covariant rule (\ref{covariant}), this amounts to consider
the cases $K=0,-1,1$.

\begin{prop}\label{dk}
Let $(M,[h_0])$ be a conformal connected manifold of dimension
$n>2$. If a conformal representative $h\in[h_0]$ has constant
sectional curvatures equal to $K=-1,0$ or $1$ then the $k$-th
conformal Laplacian $P_k$ for $h$ is
\begin{equation}\label{calcpk}
P_k=\prod_{j=1}^k\Big(\Delta_{h}+\Big(\ndemi-j\Big)\Big(\ndemi+j-1\Big)K\Big)
\end{equation}
and
\[d_k=\sharp\Big\{j\in\nn; j\leq k, -K(\ndemi-j)(\ndemi+j-1)\in\sigma(\Delta_{h})\Big\}.\]
\end{prop}
\textsl{Proof}: The flat case is well known (see \cite{GRZ}).
Since the local expression of $P_k$ on a manifold with constant
curvature $K=1$ is clearly the same as for the sphere $S^n$, the
result is obtained for example by calculating the residue of the
scattering operator at $\ndemi+k$ on the hyperbolic space, which
is $c_kP_k$ on the sphere according to \cite[Th.1]{GRZ}. Using for
example \cite[Appendix]{GZ0}, the scattering matrix on $\hh^{n+1}$
is
\[S(s)=2^{n-2\la}\frac{\Gamma(\ndemi-s)}{\Gamma(s-\ndemi)}
\frac{\Gamma\left(\sqrt{\Delta_{S^n}+(\frac{n-1}{2})^2}+\frac{1-n}{2}+s\right)}
{\Gamma\left(\sqrt{\Delta_{S^n}+(\frac{n-1}{2})^2}+\frac{n+1}{2}-s\right)}\]
and the residue is easily seen to be $c_kP_k$ with $P_k$ in
(\ref{calcpk}). To deal with the case of negative constant
curvature, one can use the expression of the scattering operator
on the cylindrical manifold studied in \cite[Appendix B]{PP} and
we find the same result. The formula for $d_k$ is a
straightforward consequence of the expression of $P_k$.
\qed\\
As a consequence, $d_k=1$ if $K=0$, $d_k\leq d_{n/2}$ for all $k$
if $K<0$ and $d_k=0$ if $k<\ndemi$ and $K>0$. The expression of
$P_k$ on the sphere was obtained previously by Thomas Branson
\cite{BR} (see also Beckner \cite{BE} when $k=\ndemi$). In a
recent preprint, Rod Gover \cite{GO} extends Proposition
\ref{calcpk} to the case
of a conformal class which contains an Einstein representative.
Note that this Proposition and consequences will not be used for what follows but 
they might be of independent interest.\\

Let us consider the counting function for resonances $\mc{N}(R)$
and for scattering poles $\mc{N}_s(R)$ defined by
$$\mc{N}(R):=\sum_{|s|\leq R,\ s\in \mc{R}}m_s, \quad
\mc{N}_s(R):=\sum_{|s|\leq R,\ s\in \mc{S}}\nu_s.$$
We clearly have
\[\mc{N}_s(R)\geq \mc{N}(R).\]
Patterson-Perry \cite{PP} proved the upper bounds $\mc{N}_s(R)=O(R^{n+1})$ which implies
in view of (\ref{multi}) (note that the term $m_{n-s_0}$ in (\ref{multi}) are non-zero 
only for finite number of terms if $s_0\in\ndemi-\nn$)
\begin{lem}
\label{rest}
Using the above notations, we have as $R\rightarrow+\infty$,
$$\mc{N}(R)=O(R^{n+1}), \quad \sum_{k=1}^R d_k=O(R^{n+1}).$$
\end{lem}
Note that those results are optimal in general in the sense that for some examples
(take $\hh^{n+1}$ with $n+1$ even for the first bound and $\hh^{n+1}$
with $n+1$ odd for second bound) this is an asymptotic.

\subsection{The hyperbolic space model}
We begin by a simple calculation for the wave kernel on the hyperbolic space $\hh^{n+1}$
which shows the relation between resonances and $0$-trace.
Recall that the set of resonances $\mc{R}_0$ for the Laplacian on $\hh^{n+1}$ is given by
\[\mc{R}_0=\left\{
\begin{array}{ll}
\{-k \textrm{ with multiplicity } h_n(k); k\in \nn_0 \} & \textrm{if }n+1 \textrm{ is even}\\
\emptyset & \textrm{if } n+1 \textrm{ is odd}
\end{array}\right. \]
where $h_n(k)$ is the dimension of the space of spherical harmonics of degree $k$ on $S^{n+1}$
\[h_n(k):=(2k+n)\frac{(k+1)(k+2)\dots (k+n-1)}{n!}, \quad h_n(0)=1\]
Let us write for simplicity
\[R_0(s):=(\Delta_{\hh^{n+1}}-s(n-s))^{-1}, \quad
U_0(t):=\cos\left(t\sqrt{\Delta_{\hh^{n+1}}-\frac{n^2}{4}}\right)\]
and $R_0(s;w,w'), U_0(t;w,w')$ their respective kernel.
Then we choose a boundary defining function $x_0$ of $\hh^{n+1}$ and define the distribution on $\rr^*$
\[u_0(t):=\zerotr(U_0(t))\]
where the $0$-trace is taken with respect to $x_0$.
\begin{lem}\label{model}
On $\hh^{n+1}$, we have the following formula for $t>0$
\[
u_0(t)=\left\{
\begin{array}{ll}
\demi\left(\sum_{k=0}^\infty h_n(k)e^{-t(\ndemi+k)}\right)
 & \textrm{if } n+1 \textrm{ is even}\\
0 & \textrm{if } n+1 \textrm{ is odd}\end{array}\right.
\]
\end{lem}

\textsl{Proof}: we use the formula of the wave kernel on $\hh^{n+1}$ (see Helgason \cite{HE}).
In odd dimension the result is clear since the wave kernel vanishes on the diagonal
\begin{equation}\label{u0odd}
U_0(t;w,w)=0.
\end{equation}
In even dimension it suffices to use the formula
\begin{equation}\label{u0even}
U_0(t;w,w)=
(-\pi)^{-\frac{n+1}{2}}n!!2^{-\frac{3(n+1)}{2}}\frac{\cosh\frac{t}{2}}{(\sinh\frac{t}{2})^{n+1}}
\end{equation}
and check that
\begin{equation}\label{formule}
\frac{\cosh\frac{t}{2}}{(\sinh\frac{t}{2})^{n+1}}=
2^n\left(\sum_{k=0}^\infty h_n(k)e^{-t(\ndemi+k)}\right)
\end{equation}
then
\[u_0(t)=\left\{
\begin{array}{ll}
(-\pi)^{-\frac{n+1}{2}}n!!2^{-\frac{n+3}{2}}\left(\sum_{k=0}^\infty h_n(k)e^{-t(\ndemi+k)}\right)
\zerov(\hh^{n+1}) & \textrm{if } n+1 \textrm{ is even}\\
0 & \textrm{if } n+1 \textrm{ is odd}\end{array}\right.
\]
and using that $\zerov(\hh^{n+1})=\frac{(-2\pi)^{\frac{n+1}{2}}}{n!!}$ when $n+1$ is even (see \cite{GR})
this gives the result.
\qed\\

\subsection{Proof of Theorem \ref{poisson}}
We will now study the general case of a convex co-compact quotient of $\hh^{n+1}$.
Let $X=\Gamma\backslash \hh^{n+1}$ be the quotient of the hyperbolic space by
a convex co-compact group of isometries.
We first define the distribution on $\rr^*$
\[u_c(t):=\zerotr\left(P_X\cos\left(t\sqrt{\Delta_{X}-\frac{n^2}{4}}\right)\right)\]
with $P_X$ the projector on the continuous part of $\Delta_X$.
Note that in Theorem \ref{poisson}, $u(t)$ is the sum of $u_c(t)$ with the contribution of the discrete spectrum, 
that is \[\sum_{\substack{s\in\mc{R}\\
\Re(s)>\ndemi}}
m_s\cosh(t(\ndemi-s)).\]
From \cite[eq. 4.13]{GZ1}\footnote{We believe there is a sign typo}, we have for $\varphi\in C_0^\infty(\rr^*)$
\begin{equation}\label{fourier}
\int u_c(t)\varphi(t)dt = (4\pi)^{-1}\int_\rr \hat{\varphi}(z)\theta(z)dz
\end{equation}
where
\[\theta(z):=
2iz\left(\zerotr\left(R_X\left(\ndemi+iz\right)-R_X\left(\ndemi-iz\right)\right)\right)\]
will be shown to be a tempered distribution on $\rr$. We recall that
the Selberg zeta function $Z(s)$ of $X=\Gamma \backslash \hh^{n+1}$ is
defined by the Euler product (which converges for $\Re(s)>\delta$),
$$Z(s)=\prod_{\gamma \in \mc{P}} \prod_{k_1,\ldots,k_n \in \nn}
\left(1-\alpha_1(\gamma)^{k_1}\ldots\alpha_n(\gamma)^{k_n}e^{-(s+|k|)l(\gamma)}\right),$$
where $|k|=k_1+\ldots+k_n$.

Let $F$ be a fundamental domain of $\Gamma$ in $\hh^{n+1}$.
Then from Patterson-Perry formula \cite[eq. 6.7]{PP} we have
\begin{equation}\label{perrypat}
\theta(z) = \theta_1(z)+\theta_2(z)
\end{equation}
\[\theta_1(z):=\frac{Z'(\ndemi+iz)}{Z(\ndemi+iz)}+\frac{Z'(\ndemi-iz)}{Z(\ndemi-iz)}\]
\[\theta_2(z):= FP_{\eps\to 0}\int_{\{x>\eps\}\cap F}2iz\left[R_0\left(\ndemi+iz;w,w'\right)-
R_0\left(\ndemi-iz;w,w'\right)\right]_{w=w'}dvol(w)\]
Using the factorization of the Selberg zeta function in a Hadamard product by Patterson-Perry \cite[Th. 1.9]{PP}
(they also use Bunke-Olbrich results \cite{BO}), we observe that
\begin{eqnarray*}
\theta_1(z)&=&-i\left\{P(z)-\chi(X)
\sum_{k=0}^\infty h_n(k)\left(\frac{1}{z-i(\ndemi+k)}-\frac{1}{z+i(\ndemi+k)}+Q(z,-k)\right)\right.\\
 & &\left. +\sum_{s\in\mc{S}}\nu_s\left(\frac{1}{z-i(\ndemi-s)}-\frac{1}{z+i(\ndemi-s)}+Q(z,s)\right)\right\}
\end{eqnarray*}
with $P(z),Q(z,s)$ some polynomials in $z$
of degree less or equal to $n$ and $\chi$ the Euler characteristic of $X$. This proves, using the arguments of
Lemma 4.7 of \cite{GZ}, that $\theta_1$ is a tempered distribution on $\rr$, hence we can define its
Fourier transform $\hat{\theta_1}$ which is also a tempered distribution on $\rr$.
Now we differentiate the last equation $n+1$ times and obtain
\begin{eqnarray*}
(-i\pl_z)^{n+1}\theta_1(z)&=&i^n(n+1)!\left\{-\chi(X)
\sum_{k=0}^\infty h_n(k)\left(\frac{1}{(z-i(\ndemi+k))^{n+2}}-\frac{1}{(z+i(\ndemi+k))^{n+2}}\right)
\right.\\
& &\left.
+\sum_{s\in\mc{S}}\nu_s\left(\frac{1}{(z-i(\ndemi-s))^{n+2}}-\frac{1}{(z+i(\ndemi-s))^{n+2}}\right)\right\}\\
\end{eqnarray*}
We combine this formula with the following
\[\mc{F}^{-1}_{t\to z}(t^{n+1}e^{i\zeta |t|})=(2\pi)^{-1}(n+1)!i^n\left(\frac{1}{(z-\zeta)^{n+2}}
-\frac{1}{(z+\zeta)^{n+2}}\right)\]
for $\Im(\zeta)\geq 0$ to conclude that
\begin{equation}\label{formule1}
t^{n+1}\hat{\theta_1}(t)=-2\pi t^{n+1}\left\{\chi(X)\left(\sum_{k=0}^\infty h_n(k)e^{-(\ndemi+k)|t|}\right)
-\sum_{s\in\mc{S}}\nu_s e^{-(\ndemi-s)|t|}\right\}
\end{equation}
in view of the identity
\[t^{n+1}\hat{\theta_1}(t)=\mc{F}_{z\to t}\left((-i\pl_z)^{n+1}\theta_1(z)\right).\]

To study $\theta_2$ we use that (see \cite{PP})
\[\theta_2(z)=\pi^{-\ndemi}\frac{\Gamma(\ndemi)}{\Gamma(n)}\frac{\Gamma(\ndemi+iz)\Gamma(\ndemi-iz)}{\Gamma(iz)\Gamma(-iz)}
\zerov(X)\]
which is a tempered distribution for $z\in\rr$ since it is bounded by a polynomial (note that the $0$-volume of $X$
could depend of the defining function $x$ when $n+1$ is odd).
We deduce that $\theta$ is a tempered distribution on $\rr$ and by (\ref{fourier}), its Fourier
transform is a tempered distribution on $\rr$ which, when restricted to $\rr^*$, is
$4\pi u_c$.
To find the Fourier transform of $\theta_2$, we remark by using again \cite[eq. 4.13]{GZ}
that this is equivalent to calculate
\[FP_{\eps\to 0}\int_{F\cap \{x>\eps\}} U_0(t;w,w)dvol(w)=4\pi \zerov(X)\left\{
\begin{array}{ll}
0 & \textrm{if } n+1 \textrm{ is odd}\\
\frac{n!!2^{-\frac{3(n+1)}{2}}}{(-\pi)^{\frac{n+1}{2}}}\frac{\cosh\frac{t}{2}}{(\sinh\frac{t}{2})^{n+1}}
&  \textrm{if } n+1 \textrm{ is even}
\end{array}\right.
\]
where we used (\ref{u0odd}) and (\ref{u0even}).
But from Epstein formula for the $0$-volume of $X$ in \cite{PP}
\[\zerov(X)=\frac{(-2\pi)^{\frac{n+1}{2}}}{n!!}\chi(X),\]
and (\ref{formule}) we deduce that
\begin{equation}\label{theta2hat}
\hat{\theta_2}(t)=\left\{
\begin{array}{ll}
0 & \textrm{if }n+1 \textrm{ is odd}\\
2\pi \chi(X)\left(\sum_{k=0}^\infty h_n(k)e^{-(\ndemi+k)|t|}\right) & \textrm{if }n+1 \textrm{ is even}
\end{array}\right.
\end{equation}
To achieve the proof of the formula relating $0$-trace and resonances,
it suffices to add (\ref{formule1}) with (\ref{theta2hat}), use
(\ref{fourier}) and we obtain, as distribution on
$\rr^*$,
\[u(t)=
\demi\sum_{s\in\mc{S}}\nu_s e^{-(\ndemi-s)|t|}+\sum_{\substack{s\in\mc{R}\\
\Re(s)>\ndemi}}
m_s\cosh(t(\ndemi-s))\]
if $n+1$ is even and
\[u(t)=
\demi\sum_{s\in\mc{S}}\nu_s e^{-(\ndemi-s)|t|}-
2^{-n-1}\chi(X)\frac{\cosh\frac{t}{2}}
{(\sinh\frac{|t|}{2})^{n+1}}
+\sum_{\substack{s\in\mc{R}\\
\Re(s)>\ndemi}}
m_s\cosh(t(\ndemi-s))\]
if $n+1$ is odd. Now it suffices to use the formula (\ref{multi}) relating $\nu_{s_0}$
and $m_{s_0}, m_{n-s_0}$. Note that
we do not analyze the singularity at $t=0$ (involving
distributions supported at $t=0$), this would require better results about
the factorization of $Z(s)$ as a product over its zeros. 
The analysis of $u(t)$ at $t=0$ is done in a general setting
by Joshi-Sa Barreto \cite[Prop. 4.3]{JSB1} and one could get actually good informations
of this singularity in our case if $\delta<n/2$.\\

To obtain the part with the length spectrum, we remark that for
$\Re(s)>\delta$,
\[Z(s)=\exp\left(-\sum_{\gamma}\sum_{m=1}^{\infty}\frac{1}{m}\frac{e^{-s ml(\gamma)}}{G_\gamma(m)}\right).\]

As before $s=\ndemi+iz$, so when the sum converges (this is the case if $\Im(z)=0<\ndemi-\delta$),
\[\theta_1(z)=\sum_{\gamma}\sum_{m=1}^{\infty}l(\gamma)G_\gamma(m)^{-1}(e^{-(\ndemi+iz)ml(\gamma)}+
e^{-(\ndemi-iz)ml(\gamma)})\]
thus
\[\hat{\theta_1}(t)=\sum_{\gamma}\sum_{m=1}^{\infty}2\pi l(\gamma)G_\gamma(m)^{-1}e^{-\ndemi ml(\gamma)}
(\delta(t+ml(\gamma))+\delta(t-ml(\gamma)))\]
and we are done at least when $\delta<\ndemi$ since there are no $L^2$-eigenvalues in this case.
The other cases are treated by the arguments of Perry \cite[lemma 2.3]{P}
using a contour deformation which makes the term
\[-\sum_{\substack{s\in\mc{R}\\ \Re(s)>\ndemi}}m_s \cosh(t(\ndemi-s))\]
appear in addition and this term cancels the term $\zerotr((1-P_X)\cos(t\sqrt{\Delta_X-\frac{n^2}{4}}))$.
\qed\\

\section{Asymptotics of the counting function for prime geodesics}
In this section, we prove Theorem \ref{lsp}. The proof is based
directly on the trace formula. The standard method using zeta functions
and contour deformation could be applied, but a good estimate of the
growth of $|Z'(s)/Z(s)|$ in strips parallel to the imaginary axis is
lacking and the error term obtained might be marginally large.
In the special case of Schottky manifolds, Guillop\'e, Lin and Zworski
\cite{glz} have shown that $|Z(s)|=O(e^{C|\Im(s)|^\delta})$. In that direction,
it could be interesting in that case to use the method of regularization, contour
deformation and then finite differences and compare the error term
obtained with the remainder of Theorem \ref{lsp}.

\bigskip
Let us define for $x\geq 1$, the following counting
functions.
$$\Pi_0(x):=\#\{\gamma \in \mc{P}\ :\ e^{l(\gamma)}\leq x \},\quad
\Pi(x):=\sum_{\substack{\gamma\in\mc{P},k\in\nn,\\ e^{kl(\gamma)}\leq x}}\frac{1}{k},\quad
\Psi(x):=\sum_{\substack{\gamma\in\mc{P},k\in\nn,\\ e^{kl(\gamma)}\leq x}}\frac{l(\gamma)}{G_\gamma(k)}$$
We have obviously $N(T)=\Pi_0(e^T)$ and
$\Pi(x)=\sum_{k=1}^{+\infty}\frac{1}{k}\Pi_0(x^{1/k})$.
Using the asymptotic formula of Perry \cite{P4}, we have
$$\Pi_0(x)=O\left(\frac{x^\delta}{\log(x)}\right),$$
which is enough to deduce
$$\Pi(x)=\Pi_0(x)+O(x^{\delta/2}).$$
The counting function $\Psi(x)$ is related to $\Pi(x)$ by remarking
that $$\int_2^x\frac{d\Psi(u)}{\log (u)}=\Pi(x)+\Phi(x)+O(1),$$
where the remainder $\Phi(x)=\sum_{kl(\gamma)\leq
  \log(x)}\frac{1}{k}(G_\gamma(k)^{-1}-1)$ can be estimated by
$$|\Phi(x)|\leq C_X\sum_{kl(\gamma)\leq \log(x)}
\frac{1}{k}e^{-kl(\gamma)}=O\left(\int_2^x\frac{d\Pi(u)}{u} \right).$$
Since we have $\Pi(x)=O(x^\delta)$, a straightforward Stieltjes
integration by parts shows that $$\Phi(x)=O(x^{\delta-1}).$$
In a nutshell, we have
\begin{equation}
\label{lsp1}
\Pi_0(x)= \int_2^x\frac{d\Psi(u)}{\log (u)}+O\left( x^{\max(\frac{\delta}{2},\delta-1)}\right).
\end{equation}
Our goal is now to get precise asymptotics of $\Psi(x)$. For this
purpose we need to introduce the following family of test functions.
Let $l(X)$ denote the length of the shortest closed geodesic on $X$.
In the following of the proof, $x,y$ will be large real parameters
satisfying $$y>0,\ 0<l(X)<x<x+y\ \textrm{and}\ y=O(x^\alpha),$$
where $\alpha<1$ will be chosen {\it a posteriori} to optimize the
error term. Let $\varphi_{x,y}$ be a $C_0^\infty(\rr)$ positive {\it even} test
function such that
$$\varphi_{x,y}(u)=\left \{ \begin{array}{lll}
0&\textrm{if}\ u \in&\left[0,\frac{l(X)}{2}\right]\\
1&\textrm{if}\ u \in&[l(X),\log(x)]\\
0&\textrm{if}\ u \in&[\log(x+y),+\infty)
\end{array} \right.$$
Clearly such a function does exist and it can be chosen such that for
all $k\in\nn$, there exists a constant $C_k>0$ such that
\begin{equation}
\label{test}
\sup_{\log(x)\leq u \leq \log(x+y)}|\varphi^{(k)}_{x,y}(u)|\leq C_k\left(\frac{x}{y}\right)^k.
\end{equation}
Set $h_{x,y}(t):=e^{\ndemi|t|}\varphi_{x,y}(t)$. Testing the trace
formula on $h_{x,y}$, we obtain the relation
\begin{equation}
 \label{trf}
 \sum_{s\in \mc{R}}
 m_s\int_0^{\infty}e^{st}\varphi_{x,y}(t)dt+\sum_{k\in \nn}
 d_k\int_0^{+\infty} e^{-kt}h_{x,y}(t)dt=
\sum_{0\leq kl(\gamma)\leq \log(x+y)}\frac{l(\gamma)}{G_\gamma(k)}\varphi_{x,y}(kl(\gamma))
\end{equation}
$$+\left(2^{-n}A(X)+2B(X)\right)\int_0^{+\infty}\frac{\cosh(t/2)}{(\sinh(t/2))^{n+1}}h_{x,y}(t)dt.$$
Obviously, we have
$$\int_0^{+\infty}\frac{\cosh(t/2)}{(\sinh(t/2))^{n+1}}h_{x,y}(t)dt=O(\log(x)),$$
and since $t \mapsto \sum_{k\in \nn} d_ke^{-kt}$ is uniformly convergent
and bounded on $[\frac{l(X)}{2},\infty )$, we can write
$$\sum_{k\in \nn} d_k\int_0^{+\infty} e^{-kt}h_{x,y}(t)dt=O\left(\log(x) x^{\ndemi} \right).$$
Since we assume that $\delta>\ndemi$, we know that the point spectrum
of the Laplacian is a non-empty finite subset of
$(0,\frac{n^2}{4})$ whose bottom is the simple eigenvalue
$\delta(n-\delta)$. Let $$\ndemi<\alpha_0\leq \alpha_1 \leq \ldots
<\alpha_p=\delta$$ be the corresponding resonances with respect to the
modified spectral parameter $s(n-s)$. We have the partition
$$\mc{R}=\{\alpha_0,\ldots,\delta \}\cup \{s\in \mc{R}\ :\ \Re(s)\leq
\ndemi \}.$$ In the following, we will denote by $\mc{R}^+=\{s\in \mc{R}\ :\ \Re(s)\leq
\ndemi \}$. For all $s\in \mc{R}$, we set
$$\psi_{x,y}(s)=\int_{0}^{+\infty}e^{st}\varphi_{x,y}(t)dt.$$
On the spectral side of formula (\ref{trf}), we have
$$\sum_{s\in \mc{R}}m_s\psi_{x,y}(s)=\sum_{s\in \mc{R}^+}m_s\psi_{x,y}(s)+\sum_{k=0}^p\psi_{x,y}(\alpha_k).$$
For all $0\leq k\leq p$, we get directly
$$\psi_{x,y}(\alpha_k)=\frac{x^{\alpha_k}}{\alpha_k}+O\left(
\frac{y}{x}x^{\alpha_k} \right)=\frac{x^{\alpha_k}}{\alpha_k}+O\left( \frac{y}{x^{\delta-1}}\right).$$
It remains to estimate carefully the spectral sum $\sum_{s\in
  \mc{R}^+}m_s\psi_{x,y}(s)$. If $s=0$, then we certainly have
$\psi_{x,y}(0)=O\left(x^{\ndemi}\right).$ If $s\neq 0$ then
integrating by parts $k$ times yields (we use the fact that
$\Re(s)\leq \ndemi$ and the estimate (\ref{test})),
\begin{equation}
\label{ipp}
\psi_{x,y}(s)=\frac{(-1)^k}{s^k}\int_0^{+\infty}e^{st}\varphi^{(k)}_{x,y}(t)dt
=O_k\left(\frac{x^{\ndemi}}{|s|^k}\left(\frac{x}{y}\right)^{k-1} \right).
\end{equation}
\noindent Writing
\begin{equation}
 \sum_{s\in \mc{R}^+}m_s\psi_{x,y}(s)=\sum_{|s|\leq
 \frac{x}{y}}m_s\psi_{x,y}(s)+\sum_{|s|\geq \frac{x}{y}}m_s\psi_{x,y}(s),
\end{equation}
and using the estimate (\ref{ipp}) for $k=1$ and $k=n+2$, we get
$$\left| \sum_{s\in \mc{R}^+}m_s\psi_{x,y}(s)\right|\leq
x^{\ndemi}\int_1^{\frac{x}{y}}\frac{d\mc{N}(u)}{u}+
x^{\ndemi}\left(\frac{x}{y}\right)^{n+1}
\int_{\frac{x}{y}}^{+\infty}\frac{d\mc{N}(u)}{u^{n+2}}+O(x^{\ndemi}).$$
A Stieltjes integration by parts combined with Lemma \ref{rest} shows that
$$\int_1^{\frac{x}{y}}\frac{d\mc{N}(u)}{u}=O\left( \left(\frac{x}{y}\right)^n
\right),$$ and a similar argument yields
$$\int_{\frac{x}{y}}^{+\infty}\frac{d\mc{N}(u)}{u^{n+2}}=O\left(
\frac{y}{x} \right).$$
Gathering all our previous estimates, we get for all $x,y$ large (and
satisfying the previously defined a priori properties),
$$\sum_{0\leq kl(\gamma)\leq
\log(x+y)}\frac{l(\gamma)}{G_\gamma(k)}\varphi_{x,y}(kl(\gamma))=\frac{x^\delta}{\delta}+\sum_{k=0}^{p-1}
\frac{x^{\alpha_k}}{\alpha_k}+O\left(\frac{y}{x^{1-\delta}}+\frac{x^{3\ndemi}}{y^n}\right)
+O\left((\log x)x^{\ndemi}\right).$$
Subtracting the above formula from that of $x+y$ instead of $x$ and
using the positivity of the left side, we observe that
$$\sum_{x\leq e^{kl(\gamma)}\leq
  x+y}\frac{l(\gamma)}{G_\gamma(k)}\varphi_{x,y}(kl(\gamma))
\leq \sum_{i=0}^p\frac{(x+y)^{\alpha_i}-x^{\alpha_i}}{\alpha_i}+
O\left(\frac{y}{x^{1-\delta}}+\frac{x^{3\ndemi}}{y^n}\right)
+O\left((\log x)x^{\ndemi}\right),$$
$$=O\left( \frac{y}{x^{1-\delta}} \right)+O\left(\frac{y}{x^{1-\delta}}+\frac{x^{3\ndemi}}{y^n}\right)
+O\left((\log x)x^{\ndemi}\right).$$
In other words, the above standard argument shows that we can drop the terms
over $e^{kl(\gamma)}$ without changing the asymptotics. We have thus
obtained
$$\Psi(x)=\frac{x^\delta}{\delta}+\sum_{k=0}^{p-1}
\frac{x^{\alpha_k}}{\alpha_k}+O\left(\frac{y}{x^{1-\delta}}+\frac{x^{3\ndemi}}{y^n}\right)
+O\left((\log x)x^{\ndemi}\right).$$
We recall that all the preceding estimates are valid for all
$y=O(x^\alpha)$ with $\alpha<1$. A straightforward computation shows
that the choice
of $$y=n^{\frac{1}{n+1}}x^{\frac{3}{n+1}\ndemi
  +\frac{1-\delta}{n+1}}$$
minimizes the global error term which becomes
$O\left(x^{\beta_n(\delta)}\right)$ where
$\beta_n(\delta)=\frac{n}{n+1}\left(\frac{1}{2}+\delta\right)$.
Remark that because $\ndemi<\delta$, we have
indeed the exponent $\frac{3}{n+1}\ndemi +\frac{1-\delta}{n+1}<1$.
Going back to formula (\ref{lsp1}), the final asymptotics follow by a
direct Stieltjes integration by parts.
\qed \\

\end{document}